\pageno=1
 \magnification=\magstep1
 \baselineskip=15pt
 \centerline{Latin Transversals of Rectangular Arrays}
 \centerline{by}
 \centerline{S.K.Stein}
 \vskip.3in

 Let $m$ and $n$ be integers, $ 2\leq m \leq n$. An $m$ by $n$ {\it
 array}
 consists of $mn$ cells, arranged in $m$ rows and $n$ columns, and each
 cell contains exactly one symbol. A {\it transversal} of an array
 consists
 of $m$ cells, one from each row and no two from the same column. A {\it
 latin transversal} is a transversal in which no symbol appears more than
 once. We investigte  $L(m,n)$, the largest integer such that if each
 symbol in an $m$ by $n$ array appears at most $L(m,n)$ times, then the
 array must have a latin transversal. We will obtain upper and lower
 bounds
 on $L(m,n)$ and also determine $L(2,n)$ and  $L(3,n)$.
 \vskip.2in
  Note that $L(m,n)\leq L(m+1,n)$ and that $L(m,n) \leq  L(m,n+1)$.
 The function $L$ satisfies two more inequalities, stated in Theorems 1
 and
 2.
 \vskip.2in
 {\bf Theorem 1.} For $n \leq 2m - 2$, $L(m,n) \leq n-1.$
 \vskip.2in
 The proof depends on a general construction due to E. T. Parker [6],
 illustrated for the cases $(m,n) = (4,4), (4,5)$, and $(4,6)$:

 $\matrix{1&1&4&4\cr
 2&2&1&1\cr
 3&3&2&2\cr
 4&4&3&3\cr}$ \quad
 \vskip.2in

  $\matrix{1&1&1&4&4\cr
 2&2&2&1&1\cr
 3&3&3&2&2\cr
  4&4&4&3&3\cr}$ \quad
 \vskip.2in

 $\matrix{1&1&1&4&4&4 \cr
 2&2&2&1&1&1 \cr
  3&3&3&2&2&2\cr
  4&4&4&3&3&3\cr}$
 \vskip.2in
       In each case an attempt to construct a latin transversal might as
 well begin with a $1$ in the top row.  The $2$ must then be selected
 from
 the $2$'s in the second row, and the $3$ from the $3$'s in the third
 row.
 Such choices do not extend to a latin transversal.

 	{\bf Theorem 2.} For $n \geq 2m-1, L(m,n) \leq {(mn - 1)/(m-1)}.$

      This theorem follows from the fact that if only $m-1$ distinct
 symbols appear in an $m$ by $n$ array, the array cannot have a latin
 transversal. In detail, if each of $m-1$ symbols appears at most $k$
 times
 and $(m - 1)k$ is at least $mn - 1$, the symbols may fill all the cells.
 Hence the inequality stated in Theorem 2 holds.

 	Though Theorem 2 is valid for all $n$, in view of Theorem 1, it is
 of interest only for $n  \geq 2m-1.$.

 	It is quite easy to determine $L(2,n)$. A moment's thought shows
 that $L(2,2) = 1$ and that $L(2,n) = 2n-1$ for $n \geq 3$. This means
 that
 for $m=2$ the inequalities in Theorems 1 and 2 become equalities.
 	The case $m=3$ is similar, for it turns out that $L(3,n)$ equals
 $n-1$ for $n = 3,4$ and is the largest integer less than or equal to
 $(3n-1)/2$ for $n \geq 5.$  The following two lemmas are used in the
 proof of the second assertion.

 	In each case $x$ stands for $1$ or $2$.

 	{\bf Lemma 1.} Assume that in a $3$ by $n$ array, $n \geq 4$, the
 following configuration is present. If there is no latin transversal the
 cell containing $y$ must be $1$ or $2$.

\vskip.1in

 $\matrix{y &\quad & \quad & \quad & \ldots \cr
 \quad & x &\quad &\quad & \ldots \cr
  \quad & \quad & 1 & 2 & \ldots \cr}$
 \vskip.1in	

The proof is immediate.
 \vskip.2in

 {\bf Lemma 2.} Assume that in a $3$ by $n$ array, $n \geq 4$, some
 symbol
 occurs at most three times.  Then, if there is no latin transversal some
 symbol occurs at least $2n - 2$ times, hence at least $3n/2$ times.

 	{\bf Proof.} We regard two arrangements of symbols in cells as
 equivalent if one can be obtained from the other by a permutation of
 rows,
 a permutation of columns, and a relabeling. There are seven inequivalent
 configurations of cells occupied by a symbol that occurs at most three
 times.  We illustrate them by treating the case when $1$ appears twice,
 in
 one row, as in the following diagram:
 \vskip.2in
 $\matrix{1 & 1 & \quad & \quad & \quad & \ldots \cr
  b&b&2 &b&b & \ldots\cr
  a&a&c&a&a&\ldots\cr}$
 \vskip.2in
 	We may assume the $2$ occurs as indicated.  It follows that the
 cells marked $a$ are filled with $2$'s.  This implies that all the cells
 marked $b$ are also filled with $2$'s, and finally  that the cell marked
 $c$ also contains a $2$. Hence the symbol $2$ appears at least $2n$
 times.

 	In the case when a symbol occurs only once, it is not hard to show
 that some symbol appears at least $2n - 2$ times. In most of the other
 cases a symbol appears almost $3n$ times.
 \vskip.2in

 	{\bf Theorem 3.} (a) $L(3,3)=2$ and $L(3,4) = 3.$ (b) For $n \geq
 5$, $L(3,n)$ is the greatest integer less than or equal to $(3n-1)/2.$

 	{\bf Proof of (a).} To begin we show that $L(3,3) = 2$.

 	By Theorem 1, we know that $L(3,3)$ is at most 2. All that remains
 is to show that if each symbol in a $3$ by $3$ array appears at most
 twice, then the array has a latin transversal.

 	First of all, at least five different symbols must appear in the
 array, hence at least one, say $1$, must appear exactly once. Without
 loss
 of generality, we may assume that the following configuration occurs in
 the array:
 \vskip.2in
 $\matrix{1&4&\quad\cr
 \quad&2&3\cr
 \quad&3&2\cr}$
 \vskip.2in

 	To avoid the formation of a latin transversal, the two empty cells
 in the first column must contain the symbol $4$, in violation of our
 assumption that each symbol appears at most twice. Thus $L(3,3)=2.$

 	Next we show that $L(3,4) = 3$.

 	We may assume that if there is no latin transversal in a $3$ by
 $4$ array that  either all cells contain the same symbol or else the
 following configuration occurs, where $x$ stands for either $1$ or $2$:
 \vskip.2in
 $\matrix{1&\quad&\quad&x\cr
 \quad&2&\quad&\quad\cr
 \quad&\quad&1&x\cr}$
 \vskip.2in
  This breaks into two cases depending on whether the top $x$ is $1$ or
 $2$. In either case, the condition that there is no latin transversal
 forces another cell to be an $x$.
 \vskip.2in
 (1) \quad $\matrix{1&\quad&\quad&1\cr
 \quad &2&\quad&\quad\cr
  x&\quad &1 & x\cr}$.
 \vskip.2in

 (2) \quad  $\matrix{1&\quad & \quad & 2\cr
  x &2 &\quad & \quad \cr
 \quad & \quad & 1 &x\cr}$

 	In case (1), since there are already three $1$'s, the $x$'s must
 be $2$'s, and the cell labeled $y$ below is part of a latin transversal.
 \vskip.2in
 $\matrix{1 & \quad & \quad  & 1 \cr
 \quad & 2  &y & \quad \cr
 2 & \quad &1 &  2 \cr}$

 Case (2) is slightly different. Using the lower $x$, we obtain two
 cases:
 \vskip.2in
 $\matrix{1  & \quad &x &  2\cr
  x &2 & \quad &\quad\cr
  \quad & \quad & 1 & 1 \cr}$

 and
 \vskip.2in

 $\matrix{1 & \quad & \quad\ &2\cr
  x & 2 & x & \quad\cr
  \quad & \quad & 1 &2 \cr}$

      In the first case both  $x$'s must be $2$'s, forcing the presence
 of
 more than three $2$'s. In the second, both $x$'s must be $1$'s, and
 there
 are at least four $1$'s. Since no symbol is assumed to appear more than
 $3$ times, these contradictions complete the proof.
 \vskip.2in

 	{\bf Proof of (b).} We show first that $L(3,5) = (3\cdot 5 -
 1)/2$, that is, $L(3,5) = 7.$

 	Consider a $3$ by $5$ array without a latin transversal. We will
 show that some symbol occurs at least eight times.

 	First of all, if all fifteen cells contain the same symbol, a
 symbol appears at least eight times. We therefore assume that at least
 two
 different symbols occur, and therefore can assume that the following
 configuration is present:
 \vskip.2in
 $\matrix{1 & \quad & \quad & x & x\cr
  \quad & 2 &\quad & \quad & \quad \cr
 \quad & \quad & 1 & \quad & \quad \cr}$

 This breaks into three cases, depending on the top two $x$'s.
 \vskip.2in
(1) $\matrix {1 & \quad & \quad & 1 & 1  \cr
  \quad & 2 & \quad& \quad & \quad \cr
 \quad & \quad & 1 & x & x\cr}$
 \vskip.2in
 (2) $\matrix { 1 & \quad & \quad & 1 & 2\cr
  \quad & 2 & \quad & x &\quad \cr
  \quad & \quad & 1 & x & x\cr}$
 \vskip.2in
 (3) $\matrix{ 1 &\quad & \quad & 2 & 2 \cr
 x & 2 & \quad & x & x \cr
 \quad & \quad & 1 & x & x \cr}$
 \vskip.2in

 We will analyze case (2); the other two cases are similar.
 	We may fill in three more cells with $x$'s, in one case using
 Lemma 1:
 \vskip.2in
 $\matrix {1 & \quad & \quad & 1 & 2 \cr
  x & 2 & \quad &x & \quad \cr
  x & x & 1 & x &x \cr }$

 This case breaks into two cases, depending on whether the lowest cell in
 the first column is filled with a $1$ or a $2$. These two cases,
 including
 the implied $x$'s in other cells are shown here:
 \vskip.2in

 $\matrix {1 & \quad & x &1 &2 \cr
 x & 2&x &x & \quad \cr
 1 & x & 1 & x &x \cr}$
 
 \vskip.2in

 $\matrix { 1 & \quad & \quad &1 &2 \cr
 x & 2 & x & x & x \cr
 2 & x & 1 & x & x \cr}$
 \vskip.2in

 In both cases two cells remain empty. If both are filled with $1$'s or
 $2$'s, then some symbol occurs more than seven times. On the other hand,
 if some other symbol occurs there, then again, by Lemma 2, some symbol
 occurs more than seven times.

 	Cases (1) and (3) are similar.  Thus $L(3,5) = 7.$

 	We will now prove by induction that for even $n \geq 6$, $L(3,n) =
 (3n-2)/2$ and that for odd $n \geq 5$, $L(3,n)= (3n-1)/2$.

 	Assume that the induction holds for a particular even $n$, that
 is, $L(3,n)=(3n-2)/2.$  We will show that it holds for $n+1$, which is
 odd, that is, $L(3,n+1)= (3n+2)/2.$ Note that in this case we would have
 $L(3,n+1)=L(3,n)+ 2.$

 Consider a $3$ by $n+1$ array in which each symbol occurs at most
 $(3n+2)/2$ times. If each symbol occurs at most $(3n-2)/2$ times, delete
 one column, obtaining a $3$ by $n$ array, which has a latin transversal,
 by the inductive assumption. Hence the original array has, also.

 	Now assume that there is at least one symbol occurring at least
 $3n/2$ times.  If there are two such symbols, they occupy at least $3n$
 cells.  Hence some symbol appears at most three times.  By Lemma 2, some
 symbol occurs at least $3(n+1)/2$ times, which contradicts the
 assumption
 that each symbol occurs at most $(3n+2)/2$ times.

 Hence there is only one symbol that occurs at least $3n/2$ times, that
 is,
 $3n/2$ or $(3n+2)/2$ times. There must be a column in which it appears
 at
 least twice.  Deleting that column, we obtain a $3$ by $n$ array in
 which
 each symbol occurs at most $(3n-2)/2$ times.  By the inductive
 assumption,
 this array has a latin transversal, hence the original array does.

 The argument when $n$ is odd and $n+1$ is even is similar. (It is a bit
 shorter, since in this case $L(3,n+1)= L(3,n)+1.$) This completes the
 proof of the theorem.
 \vskip.1in
 	The next theorem gives a non-trivial lower bound on $L(m,n)$.
 \vskip.1in
 {\bf Theorem 4.} $L(m,n)\geq n-m+1$.
 \vskip.1in
 The proof is an induction on $m$.
 \vskip.1in
 The theorem is true for $m=2$ or $m=3$ and for $n=m.$ We will consider
the case $n\geq m+1.$  Assuming the theorem is true for $m-1$,
 we
 will prove it for any array $A$ with $m$ rows.  In order to simplify the
 diagrams and the exposition, we consider the case $m=5$, which
 illustrates
 the argument in the general case.
 \vskip.1in
 Assuming that $L(4,n)$ is at least $n-3$, we will show that $L(5,n)$ is
 at
 least $n-4$.
 \vskip.1in
 	Consider a $5$ by $n$ array $A$ in which each symbol appears at
 most $n-4$ times.  By the induction assumption, the $4$ by $n$ array
 consisting of the first four rows of $A$ has a latin transversal.
 Assuming that $A$ does not have a latin transversal, we may conclude
 that
 $A$ contains an equivalent of the following configuration:
\vskip.2in

 $\matrix{1&y&y&y&y&x&x&x&x&.&.&.\cr
 \quad&2&\quad &\quad&\quad&\quad&\quad&\quad&\quad&.&.&.\cr
 \quad&\quad&3&\quad&\quad&\quad&\quad&\quad&\quad&\quad&.&.&.\cr
 \quad&\quad&4&\quad&\quad&\quad&\quad&\quad&.&.&.\cr
 y&y&y&y&1&x&x&x&x&.&.&.\cr}$
\vskip.2in

 An $x$ stands for $1,2,3$ or $4$ while a $y$ stands for any symbol in
 $A$.
 There are cells marked $x$ since we are assuming that $A$ has no latin
 transversal.

 At this point $2(n-4)$ cells contain $x$ or $1$. Since $1$ occurs
at most $n-4$ times in $A$, there must be a $2,3,$ or $4$ in some cell
 marked
 $x$.  It is no loss of generality to take that symbol to be $2$.

 No matter which $x$ is replaced by $2$, there is a unique partial
 latin transversal consisting of that cell and cells marked $1,3,$ and
$4$.
 This
 permits us to fill in the second row with four $y$'s, one in each column
 that meets that transversal, and $n-5$ $x$'s.  The next diagram
 illustrates one of the two cases. \vskip.1in
\vskip.2in

 $\matrix{1&y&y&y&y&2&x&x&x&.&.&.\cr
 x&2&y&y&y&y&x&x&x&.&.&.\cr
 \quad&\quad&3&\quad&\quad&\quad&\quad&\quad&\quad&.&.&.\cr
 \quad&\quad&\quad&4&\quad&\quad&\quad&\quad&\quad&.&.&.\cr
 y&y&y&y&1&x&x&x&x&.&.&.\cr}$
 \vskip.2in

There are now $3(n-4)$ cells containing $x,1$, or $2$.  Among these
 cells
 must be a cell containing either $3$ or $4$.  We may assume, after
 permuting rows and columns and relabeling, that it is $3.$  There are
 essentially five different positions in which that symbol may appear,
 depending on which of the three rows that have $x$'s it lies in and
 where,
 relative to cells containing $1$ or $2$, it is situated.  (In each
case
 a
 unique partial latin transversal forms with cells labeled $1,2,3,$ or
$4$.)
 The
 following diagram illustrates one case. The $y$'s in the third row are
 again in the columns that contain the partial  latin transversal that includes
 the cell with the new $3$, which is not in the third row.
 \vskip.1in
 $\matrix{1&y&y&y&y&2&x&x&x&.&.&.\cr
 x&2&y&y&y&y&x&x&x&.&.&.\cr
 y&y&3&y&x&y&x&x&x&.&.&.\cr
 \quad&\quad&\quad&4&\quad&\quad&\quad&\quad&\quad&.&.&.\cr
 y&y&y&y&1&3&x&x&x&.&.&.\cr}$
 \vskip.1in
 There are now at least $4(n-4)+1$ cells occupied by $1,2,3,$ or
 $4$.  Since each symbol appears at most $n-4$ times, this is a
 contradiction, and the theorem is proved.

 In view of our experience with $m=2$ and $3$, it is tempting to
 conjecture
 that the values for $L(m,n)$ suggested by Theorems 1 and 2 would be
 correct even for $m \geq 4.$ In other words, one is tempted to
 conjecture
 that for $ m \leq n \leq 2m-2$, we have $L(m,n) = n-1$ and that for $n
 \geq 2m-1$, we have $L(m,n)$ equal to the greatest integer less than or
 equal to $(mn-1)/(m-1)$.
 Dean Hickerson has shown that $L(4,4)= 3$, in agreement with the first
 part of the conjecture. However, he also has shown that $L(4,7)$ is at
 most $8$, hence is a counterexample to the second part.

 	A result of Hall [4] lends some support for the conjecture that
 $L(n-1,n)=n-1$. Consider an abelian group of order $n, A = \{a_1,
 a_2,\ldots, a_n\}$ and  $b_1, b_2,\ldots, b_{n-1}$, a sequence of
 $n - 1$ elements of $A$, not necessarily distinct. Construct an
 $n - 1$ by $n$ array by placing $b_{i}a_j$ in the cell where row $i$
 meets
 column $j$. Hall proved that such an array has a latin transversal.
 \vskip.2in
 	Stein [9] showed that in an $n$ by $n$ array where each element
 appears exactly $n$ times there is a transversal with at least
 approximately $(0.63)n$ distinct elements.
 	Erd\''os and Spencer [3] showed that an $n$ by $n$ array in which
 each symbol appears at most $(n-1)/16$ times has a latin transversal.
 This
 gives a lower bound for $L(n,n)$, namely $(n-1)/16$. The algorithm in
 which you try to construct a latin transversal by choosing a cell in the
 top row, then a cell in the second row, and work down row by row yields
 a
 different lower bound. If each symbol appears at most $k$ times in an
 $m$
 by $n$ array, the algorithm is certainly successful if $(m-1)k \leq n -
 1.$ This implies $L(m,n) \geq (n-1)/(m-1)$.
 \vskip.2in
 	Snevily [8] offered a conjecture closely related to Hall's
 theorem: Any $k$ by $k$ submatrix of the group table of an abelian group
 of odd order has a latin transversal. (Note that in such a matrix each
 symbol appears at most $k$ times.)
 \vskip.2in
 	In the case of latin squares there are several results  concerning
 transversals that have many distinct elements, cited in [2, 9]. Ryser
 [7]
 conjectured that every latin square of odd order has a latin
 transversal,
 and, more generally, that the number of latin transversals of a latin
 square has the same parity as the order of the square. However, E. T.
 Parker [6] pointed out that many latin squares of order $7$ have an even
 number of latin transversals, for instance (6) and many other cases in
 [5].  Confirming half of Ryser's conjecture, Balasubramanian [1] proved
 that a latin square of even order has an even number of latin
 transversals.
 \vskip.2in

 {\bf References}
 \vskip.1in
 1. K. Balasubramanian, On transversals of latin squares, {\it Linear
 Algebra and Its Applications} {\bf 131} (1990) 125-129.

 2. P. Erdos, D. R. Hickerson, D.A.Norton, S. K. Stein, Has every
 latin
 square of order $n$ a partial latin transversal of size $n-1$?, {\it
 Amer.
 Math. Monthly} {\bf 95} (1988) 428-430.

 3. P. Erdos and J. Spencer, Lopsided Lovasz Local Lemma and latin
 transversals, {\it Discrete Applied Math.} {\bf 30} (1991) 151-154.

 4. M. Hall, Jr., A combinatorial problem on groups, {\it Proc. Amer.
 Math.
 Soc.} {\bf 3} (1952) 584-587.

 5. H. W. Norton, The $7$ x $7$ squares, {\it Annals of Eugenics} {\bf 9}
 (1939) 269-307.

 6. E. T. Parker, correspondence.

 7. H. J. Ryser, Neuere Problem in der Kombinatorik, in Vortraheuber
 Kombinatorik, Oberwohlfach,1967, 69-91.
 \vskip.2in
 8. H. S. Snevily, The Cayley addition table of $Z_n$, {\it Amer. Math.
 Monthly} {\bf 106} (1999) 584-585.

 9. S. Stein, Transversals of latin squares and their generalizations,
 {\it
 Pacific J. Math.} 59 (1975) 567-575.

\vskip.4in

 Mathematics Department \vskip.02in
 University of California at Davis \vskip.02in
 1 Shields Ave.\vskip.02in
 Davis, CA 95616-8633\vskip.02in
 stein@math.ucdavis.edu

 \vfill
 \end